\documentclass{llncs}
\usepackage{arxiv}

\usepackage{amsmath}
\usepackage[math]{blindtext}

\usepackage{mathtools}

\usepackage{microtype}
\usepackage[rm={oldstyle=false,proportional=true},sf={oldstyle=false,proportional=true},tt={oldstyle=false,proportional=true,variable=true},qt=false%
]{cfr-lm}
\usepackage{fancyhdr}
\usepackage{amsfonts}
\usepackage{mathtools,amssymb}
\usepackage{amscd}
\usepackage{latexsym}
\usepackage{graphicx}
\usepackage{graphics}
\usepackage{afterpage}
\usepackage[english]{babel}
\usepackage{xspace}
\usepackage{xcolor}
\usepackage{booktabs}
\usepackage{array,tabularx}
\usepackage[tableposition=top]{caption}
\usepackage{bbm}
\usepackage[]{algorithm2e}
\usepackage[
bookmarks=false,
breaklinks=true,
colorlinks=true,
linkcolor=black,
citecolor=black,
urlcolor=black,
pdfpagelayout=SinglePage,
pdfstartview=Fit
]{hyperref}

\usepackage[all]{hypcap}
\usepackage{paralist}
\usepackage{csquotes}

\usepackage{url}
\makeatletter
\g@addto@macro{\UrlBreaks}{\UrlOrds}
\makeatother

\usepackage[all]{hypcap}

\usepackage{pdfcomment}

\usepackage[capitalise,nameinlink]{cleveref}

\crefformat{footnote}{#2\footnotemark[#1]#3}
\crefname{section}{Sect.}{Sect.}
\Crefname{section}{Section}{Sections}

\usepackage{xspace}

\DeclareFontFamily{U}{MnSymbolC}{}
\DeclareSymbolFont{MnSyC}{U}{MnSymbolC}{m}{n}
\DeclareFontShape{U}{MnSymbolC}{m}{n}{
    <-6>  MnSymbolC5
   <6-7>  MnSymbolC6
   <7-8>  MnSymbolC7
   <8-9>  MnSymbolC8
   <9-10> MnSymbolC9
  <10-12> MnSymbolC10
  <12->   MnSymbolC12%
}{}
\DeclareMathSymbol{\powerset}{\mathord}{MnSyC}{180}

\usepackage{float}
\floatstyle{plaintop}
\restylefloat{table}
\restylefloat{figure}

\usepackage{tikz}
\usepackage{pgf}
\usetikzlibrary{fit,matrix,arrows,petri,topaths,positioning,automata,shapes,calc,backgrounds,chains,scopes,graphs}
\makeatletter
\tikzset{circle split part fill/.style  args={#1,#2}{%
 alias=tmp@name, 
  postaction={%
    insert path={
     \pgfextra{%
     \pgfpointdiff{\pgfpointanchor{\pgf@node@name}{center}}%
                  {\pgfpointanchor{\pgf@node@name}{east}}%
     \pgfmathsetmacro\insiderad{\pgf@x}
      \fill[#1] (\pgf@node@name.base) ([xshift=-\pgflinewidth]\pgf@node@name.east) arc
                          (0:180:\insiderad-\pgflinewidth)--cycle;
      \fill[#2] (\pgf@node@name.base) ([xshift=\pgflinewidth]\pgf@node@name.west)  arc
                           (180:360:\insiderad-\pgflinewidth)--cycle;            
         }}}}}  
 \makeatother  
 
 \makeatletter
\def\moverlay{\mathpalette\mov@rlay}
\def\mov@rlay#1#2{\leavevmode\vtop{%
   \baselineskip\z@skip \lineskiplimit-\maxdimen
   \ialign{\hfil$\m@th#1##$\hfil\cr#2\crcr}}}
\newcommand{\charfusion}[3][\mathord]{
    #1{\ifx#1\mathop\vphantom{#2}\fi
        \mathpalette\mov@rlay{#2\cr#3}
      }
    \ifx#1\mathop\expandafter\displaylimits\fi}
\makeatother

\DeclareMathOperator{\dede}{d}
\DeclareMathOperator{\dedebig}{D}
\DeclareMathOperator{\conga}{\cong_a}
\DeclareMathOperator{\Int}{Int}

\DeclareMathOperator{\I}{I}

\DeclareMathOperator{\Tr}{Tr}

\newcolumntype{C}{>{\centering\arraybackslash}X}

\begin{document}
\input glyphtounicode.tex
\pdfgentounicode=1

\title{A computation of the ninth Dedekind Number}
\author{Christian J\"akel\\[1mm]	
Technische Universit\"at Dresden, 01062, Dresden, Germany\\[1mm]
\email{christian.jaekel(at)tu-dresden.de}\\
}
\date{\vspace{-5ex}}			
\maketitle
\begin{abstract} 
In this article, we present an algorithm to compute the 9th Dedekind Number. The key aspects are the use of matrix multiplication and symmetries in the free distributive lattice, which are detected with techniques from Formal Concept Analysis.
\end{abstract}

\section{Introduction}\label{Intro}
\phantom{MI} The \emph{Dedekind numbers} are a fast growing sequence of integers, named after Richard Dedekind. Their determination is known as \emph{Dedekind's Problem}. Let $\dede(n)$ denote the $n$-th Dedekind number. A survey over past achievements can be found in \cite{YusunThesis2008}, which also contains the following table: 
\begin{figure}[H]
 \begin{center}
 	\begin{tabular}{ |c|r|c| } 
 		\hline
 		$n$ & \multicolumn{1}{|c|}{$\dede(n)$} & Source \\ 
 		\hline
 		$0$ & $2$\phantom{Mi} &  \\ 
 		$1$ & $3$\phantom{Mi} &  \\ 
 		$2$ & $6$\phantom{Mi} & Dedekind, 1897 \cite{Dedekind1897} \\
 		$3$ & $20$\phantom{Mi} &  \\
 		$4$ & $168$\phantom{Mi} &  \\
 		\hline
 		$5$ & $7581$\phantom{Mi} & Church, 1940 \cite{Church40} \\
 		\hline
 		$6$ & $7828354$\phantom{Mi} & Ward, 1946 \cite{Ward46} \\
 		\hline
 		$7$ & $2414682040998$\phantom{Mi} & Church, 1965 \cite{Church65} \\
 		\hline
 		$8$ & \phantom{Mi}$56130437228687557907788$\phantom{Mi}&\phantom{Mi} Wiedemann, 1991 \cite{Wiedemann1991} \phantom{Mi}\\
 		\hline
 	\end{tabular}
 \caption{Dedekind Numbers up to $n=8$. Integer sequence A000372.}
 \end{center}
\end{figure}
Several interpretations of $\dede(n)$ exist. For example, the value $\dede(n)$ is equal to the number of antichains of the powerset lattice $\mathbbm{2}^{n}:=(2^{\{0,\dots,n-1\}},\subseteq)$, or the number of monotone Boolean functions on $n$ variables, or the number of elements of the free distributive lattice with $n$ generators.\\

Therefore, the next section explains how the free distributive lattice can be numerically represented. After that, Section \ref{DefiNotationSection} treats (anti) isomorphic lattice intervals, a tool that is utilized for the computation. Section \ref{WarmupSection} and \ref{P4Section} deal with theoretical, and Section \ref{NinthDedekindAlgo} with practical aspects of our algorithm to compute the ninth Dedekind number.  
\section{The Free Distributive Lattice}\label{FDLSection}
Let $\mathbb{D}_n=(\dedebig(n),\leq)$ denote the free distributive lattice with $n$ generators. The following theorem can be found in \cite{FIDYTEK2001203}:
\begin{theorem}\label{LiftingTheorem}
	There is a one to one correspondence between elements of $\mathbb{D}_n$ and monotone mappings from $\mathbbm{2}^{k}$ into $\mathbb{D}_{n-k}$.
\end{theorem}
\begin{corollary}\label{Generation}
	There is a one to one correspondence of elements from $\mathbb{D}_{n+1}$ and pairs $(x,y)$ of elements $x,y$ from $\mathbb{D}_{n}$, such that $x\leq y$.
	\begin{proof}
		\emph{
		Let $x$ and $y$ be elements of $\mathbb{D}_n$. We have $\emptyset$ and $\{0\}$ as elements of $\mathbbm{2}$. Mapping $\emptyset$ to $x$ allows to monotonously map $\{0\}$ to every $y$ for which $x\leq y$ holds. \qed
		}
	\end{proof}
\end{corollary}

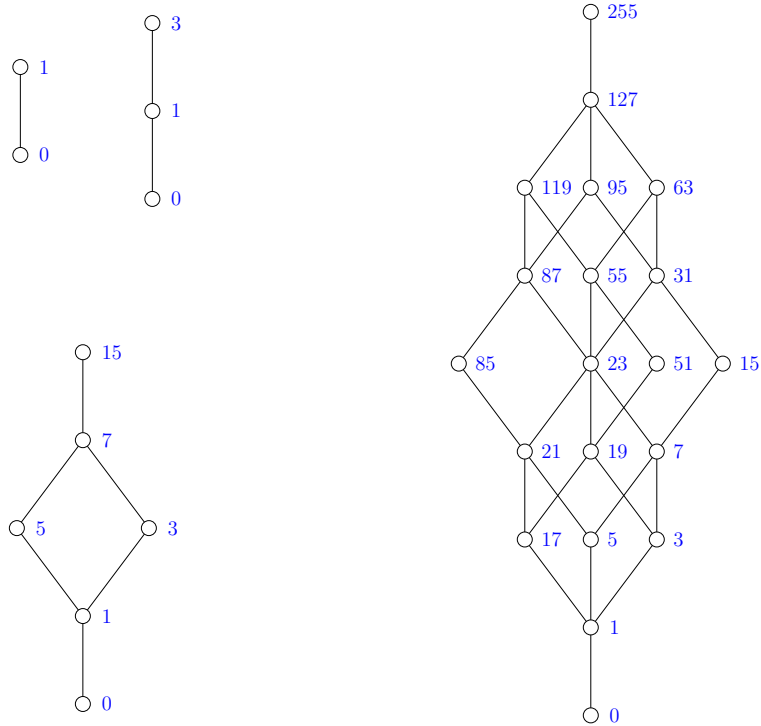
\begin{figure}[h]
	\begin{center}
		\scalebox{.65}{
			\begin{minipage}{0.49\textwidth}		
				\begin{tikzpicture}[scale=0.9,shorten >=0pt,shorten <=0pt,auto,cross line/.style={preaction={draw=white, -,line width=6pt}}]
					
					\tikzstyle{empty}=[minimum width=2pt, shape=circle , draw=black,line width=0mm]
					
					\node[empty,label={[blue]right:{ \large $\;0$}}] (00) at (0,1){\nodepart{lower}};
					\node[empty,label={[blue]right:{ \large $\;1$}}] (02) at (0,3) {\nodepart{lower}};
					\path[-](00) edge node[auto]{}(02);
					
					\node[empty,label={[blue]right:{ \large $\;0$}}] (30) at (3,0){\nodepart{lower}};
					\node[empty,label={[blue]right:{ \large $\;1$}}] (32) at (3,2) {\nodepart{lower}};
					\node[empty,label={[blue]right:{ \large $\;3$}}] (34) at (3,4) {\nodepart{lower}};
					\path[-](30) edge node[auto]{}(32);
					\path[-](32) edge node[auto]{}(34);
					
				\end{tikzpicture}\\
				\begin{minipage}{0.5\textwidth} 
					\begin{tikzpicture}[scale=0.9,shorten >=0pt,shorten <=0pt,auto,cross line/.style={preaction={draw=white, -,line width=6pt}}]
						
						\tikzstyle{empty}=[minimum width=2pt, shape=circle , draw=black,line width=0mm]
						
						\node[shape=circle , draw=white] (011) at (0,11){\nodepart{lower}};
						
						\node[empty,label={[blue]right:{ \large $\;0$}}] (00) at (0,0){\nodepart{lower}};
						\node[empty,label={[blue]right:{ \large $\;1$}}] (02) at (0,2) {\nodepart{lower}};
						\node[empty,label={[blue]right:{ \large $\;5$}}] (M24) at (-1.5,4) {\nodepart{lower}};
						\node[empty,label={[blue]right:{ \large $\;3$}}] (24) at (1.5,4) {\nodepart{lower}};
						\node[empty,label={[blue]right:{ \large $\;7$}}] (06) at (0,6) {\nodepart{lower}};
						\node[empty,label={[blue]right:{ \large $\;15$}}] (08) at (0,8) {\nodepart{lower}};
						
						\path[-](00) edge node[auto]{}(02);
						\path[-](02) edge node[auto]{}(M24);
						\path[-](02) edge node[auto]{}(24);
						\path[-](M24) edge node[auto]{}(06);
						\path[-](24) edge node[auto]{}(06);
						\path[-](06) edge node[auto]{}(08);
					\end{tikzpicture}
				\end{minipage}	
				
			\end{minipage}	
			\begin{minipage}{0.49\textwidth}
				\begin{center}
					\begin{tikzpicture}[scale=0.9,shorten >=0pt,shorten <=0pt,auto,cross line/.style={preaction={draw=white, -,line width=6pt}}]
						\tikzstyle{empty}=[minimum width=2pt, shape=circle , draw=black,line width=0mm]
						
						\node[empty,label={[blue]right:{ \large $\;0$}}] (00) at (0,0){\nodepart{lower}};
						\node[empty,label={[blue]right:{ \large $\;1$}}] (02) at (0,2) {\nodepart{lower}};
						
						\node[empty,label={[blue]right:{ \large $\,5$}}] (04) at (0,4) {\nodepart{lower}};
						\node[empty,label={[blue]right:{ \large $\,17$}}] (M154) at (-1.5,4) {\nodepart{lower}};
						\node[empty,label={[blue]right:{ \large $\,3$}}] (154) at (1.5,4) {\nodepart{lower}};
						
						\node[empty,label={[blue]right:{ \large $\,19$}}] (06) at (0,6) {\nodepart{lower}};
						\node[empty,label={[blue]right:{ \large $\,21$}}] (M156) at (-1.5,6) {\nodepart{lower}};
						\node[empty,label={[blue]right:{ \large $\,7$}}] (156) at (1.5,6) {\nodepart{lower}};
						
						\node[empty,label={[blue]right:{ \large $\,23$}}] (08) at (0,8) {\nodepart{lower}};
						\node[empty,label={[blue]right:{ \large $\,85$}}] (M38) at (-3,8) {\nodepart{lower}};
						\node[empty,label={[blue]right:{ \large $\,51$}}] (158) at (1.5,8) {\nodepart{lower}};
						\node[empty,label={[blue]right:{ \large $\,15$}}] (38) at (3,8) {\nodepart{lower}};
						
						\node[empty,label={[blue]right:{ \large $\,55$}}] (010) at (0,10) {\nodepart{lower}};
						\node[empty,label={[blue]right:{ \large $\,87$}}] (M1510) at (-1.5,10) {\nodepart{lower}};
						\node[empty,label={[blue]right:{ \large $\,31$}}] (1510) at (1.5,10) {\nodepart{lower}};
						
						\node[empty,label={[blue]right:{ \large $\,95$}}] (012) at (0,12) {\nodepart{lower}};
						\node[empty,label={[blue]right:{ \large $\,119$}}] (M1512) at (-1.5,12) {\nodepart{lower}};
						\node[empty,label={[blue]right:{ \large $\,63$}}] (1512) at (1.5,12) {\nodepart{lower}};
						
						\node[empty,label={[blue]right:{ \large $\,127$}}] (014) at (0,14){\nodepart{lower}};
						\node[empty,label={[blue]right:{ \large $\,255$}}] (016) at (0,16) {\nodepart{lower}};
						
						\path[-](00) edge node[auto]{}(02);
						\path[-](02) edge node[auto]{}(04);
						\path[-](02) edge node[auto]{}(M154);
						\path[-](02) edge node[auto]{}(154);
						
						\path[-](04) edge node[auto]{}(M156);
						\path[-](04) edge node[auto]{}(156);
						\path[-](M154) edge node[auto]{}(M156);
						\path[-](M154) edge node[auto]{}(06);
						\path[-](154) edge node[auto]{}(156);
						\path[-](154) edge node[auto]{}(06);
						
						\path[-](06) edge node[auto]{}(08);
						\path[-](M156) edge node[auto]{}(08);
						\path[-](156) edge node[auto]{}(08);
						
						\path[-](06) edge node[auto]{}(158);
						\path[-](156) edge node[auto]{}(38);
						\path[-](M156) edge node[auto]{}(M38);
						
						\path[-](010) edge node[auto]{}(08);
						\path[-](M1510) edge node[auto]{}(08);
						\path[-](1510) edge node[auto]{}(08);
						
						\path[-](010) edge node[auto]{}(158);
						\path[-](1510) edge node[auto]{}(38);
						\path[-](M1510) edge node[auto]{}(M38);
						
						\path[-](012) edge node[auto]{}(M1510);
						\path[-](012) edge node[auto]{}(1510);
						\path[-](M1512) edge node[auto]{}(M1510);
						\path[-](M1512) edge node[auto]{}(010);
						\path[-](1512) edge node[auto]{}(1510);
						\path[-](1512) edge node[auto]{}(010);
						
						\path[-](016) edge node[auto]{}(014);
						\path[-](014) edge node[auto]{}(012);
						\path[-](014) edge node[auto]{}(M1512);
						\path[-](014) edge node[auto]{}(1512);
					\end{tikzpicture}	
				\end{center}
			\end{minipage}
			\caption{The free distributive lattice with $0,1,2$ and $3$ generators.}\label{FDLFigure}
		}
	\end{center}
\end{figure}

Corollary \ref{Generation} provides a method to generate elements of $\mathbb{D}_n$ recursively. Even more, this leads to a numerical representation of $\dedebig(n)$ which facilitates the efficient implementation of algorithms to compute $\dede(n)$. For that, let the two elements of $\mathbb{D}_0$ be represented by the binary numbers $0$ and $1$. The underlying order is $0\leq 1$ applied bitwise. In the next step, we form valid pairs w.r.t. Corollary \ref{Generation} and concatenate them to $2$-bit numbers, which are $00$, $01$ and $11$. This process can be iterated further.
\begin{example} This example illustrates the generation of $\dedebig(3)$, starting from $\dedebig(0)$.\\

	\begin{minipage}{0.1\textwidth} 
		\begin{center}
			$0$\\
			$1$
		\end{center}
	\end{minipage}
	\begin{minipage}{0.05\textwidth} 
		$\rightarrow$
	\end{minipage}
	\begin{minipage}{0.1\textwidth} 
		\begin{center}
			$00$\\
			$01$\\
			$11$
		\end{center}
	\end{minipage}
	\begin{minipage}{0.05\textwidth} 
		$\rightarrow$
	\end{minipage}
	\begin{minipage}{0.1\textwidth} 
		\begin{center}
			$0000$\\
			$0001$\\
			$0011$\\
			$0101$\\			
			$0111$\\
			$1111$
		\end{center}
	\end{minipage}
	\begin{minipage}{0.05\textwidth} 
		$\rightarrow$
	\end{minipage}
	\begin{minipage}{0.54\textwidth} 
		\begin{flushleft}
			$0000.0000$   $0001.0011$   $0101.0101$\\
			$0000.0001$   $0001.0101$   $0101.0111$\\
			$0000.0011$   $0001.0111$   $0101.1111$\\
			$0000.0101$   $0001.1111$   $0111.0111$\\
			$0000.0111$   $0011.0011$   $0111.1111$\\
			$0000.1111$   $0011.0111$   $1111.1111$\\
			$0001.0001$   $0011.1111$   
		\end{flushleft}
	\end{minipage}
\end{example}
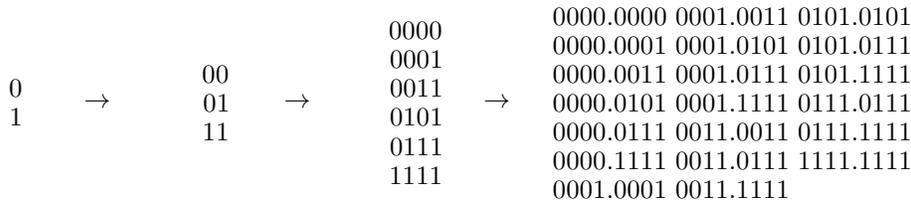
Figure \ref{FDLFigure} shows how these numbers form the free distributive lattice, by taking "bitwise or" as join- and "bitwise and" as meet operation. Instead of binary numbers, we use the more convenient decimal interpretation as integer values. From a programming point of view, this representation is memory efficient and the bitwise operations are very fast.
\section{Lattices, Definitions and Notation}\label{DefiNotationSection}
Throughout this section, we consider a finite lattice $\mathbbm{L}=(L,\vee,\wedge,\bot,\top)=(L,\leq)$ with \emph{bottom} $\bot$ and \emph{top} $\top$. An \emph{isomorphism} between lattices is a bijective map that preserves join and meet operations, and consequently order as well. Dually, an \emph{anti isomorphism} $\varphi$ is a bijective map which reverses order, and hence join and meet too. \[x\leq y \Leftrightarrow \varphi(y)\leq\varphi(x)\;\;\;\;\;\;\;\; x\vee y \Leftrightarrow \varphi(x)\wedge\varphi(y)\;\;\;\;\;\;\;\; x\wedge y \Leftrightarrow \varphi(x)\vee\varphi(y).\]
\begin{definition}
For $a,b\in L$ with $a\leq b$, we define the \emph{interval} $[a,b]:=\{x\mid x\in L,\; a\leq x\leq b\}$. An interval $I$ is \emph{isomorphic} to interval $J$, denoted by $I\cong J$, if an isomorphism exists, that maps $I$ onto $J$. Similarly, \emph{anti isomorphic} intervals $I\conga J$ are defined through the existence of a respective anti isomorphism. Two intervals $I$ and $J$ are \emph{equivalent}, denoted by $I\equiv J$, if they are isomorphic or anti isomorphic. Furthermore, $\Int(\mathbbm{L})$ denotes the set of all intervals of $\mathbbm{L}$ and $\# I$  the cardinality of $I$.
\end{definition}
\begin{proposition}
The equivalence of intervals is an equivalence relation on $\Int(\mathbbm{L})$. Consequently, the set of all intervals $\Int(\mathbbm{L})$ can be factorized by $\equiv$.
\begin{proof}
	\emph{
Reflexivity is induced by $\cong$ and the fact that every interval is isomorphic to itself. Furthermore, symmetry is given due to the symmetry of $\cong$ and $\cong_a$. To deduce transitivity, we consider three cases:
\[I\cong J \wedge J\cong K \Rightarrow I\cong K,\;\;\;\;I\conga J \wedge J\conga K \Rightarrow I\cong K\;\;\;\;I\cong J \wedge J\conga K \Rightarrow I\conga K.\]
\qed
}
\end{proof}
\end{proposition}
 For $[I]\in\Int(\mathbbm{L})/_\equiv$ and we denote the equivalence class $[I]$'s cardinality by $\#[I]$. Lastly, for every $x\in I$, we introduce two operators that decode the interval length w.r.t. the top and bottom of $I$: \[\bot_{\I}(x):=\#[\bot_{\I},x]\;\text{and}\;\top_{\I}(x):=\#[x,\top_{\I}],\]
where $\bot_{\I}$ denotes the bottom and $\top_{\I}$ the top element of $I$. Since $I$ can be deduced from the outer context (the domain containing the operators' argument), we will omit subscript $\I$ in the rest of this paper. 

\begin{example}\label{IsoIntervalExample}
These are the equivalence classes of $\Int(\mathbbm{D}_2)/_\equiv$ (compare Fig. \ref{FDLFigure}).
\begin{align*}
& \{[0,0],[1,1],[3,3],[5,5],[7,7],[15,15]\}, \{[0,1],[1,3],[1,5],[3,7],[5,7],[7,15]\},\\
& \{[0,3],[0,5],[3,15],[5,15]\}, \{[1,7]\}, \{[0,7],[1,15]\}, \{[0,15]\}.
\end{align*}  
\end{example}

\section{Enumeration of $\mathbb{D}_{n+1}$, $\mathbb{D}_{n+2}$ and $\mathbb{D}_{n+3}$ }\label{WarmupSection}
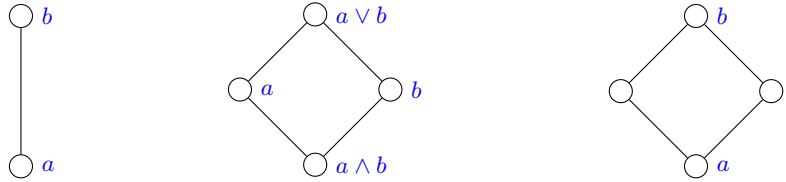
\begin{figure}
	\begin{center}
		\begin{minipage}{0.3\textwidth}
		\begin{center}		
			\begin{tikzpicture}[shorten >=0pt,shorten <=0pt,auto,cross line/.style={preaction={draw=white, -,line width=6pt}}]
				\node[shape=circle,minimum size=0.2cm , draw=black,label={[blue]right:{$a$}}] (a) at (0,0){};
				\node[shape=circle,minimum size=0.2cm , draw=black,label={[blue]right:{$b$}}] (b) at (0,2){};				
				\path[-](a) edge node[auto]{}(b);
			\end{tikzpicture}			
		\end{center}
		\end{minipage}
		\begin{minipage}{0.3\textwidth} 
			\begin{tikzpicture}[shorten >=0pt,shorten <=0pt,auto,cross line/.style={preaction={draw=white, -,line width=6pt}}]
				\node[shape=circle,minimum size=0.2cm , draw=black,label={[blue]right:{$a\wedge b$}}] (a) at (0,0){};
				\node[shape=circle,minimum size=0.2cm , draw=black,label={[blue]right:{$b$}}] (b) at (1,1){};
				\node[shape=circle,minimum size=0.2cm , draw=black,label={[blue]right:{$a$}}] (c) at (-1,1){};
				\node[shape=circle,minimum size=0.2cm , draw=black,label={[blue]right:{$a\vee b$}}] (d) at (0,2){};
				
				\path[-](a) edge node[auto]{}(b);
				\path[-](c) edge node[auto]{}(d);
				\path[-](b) edge node[auto]{}(d);
				\path[-](a) edge node[auto]{}(c);
			\end{tikzpicture}
		\end{minipage}
		\begin{minipage}{0.3\textwidth} 
			\begin{tikzpicture}[shorten >=0pt,shorten <=0pt,auto,cross line/.style={preaction={draw=white, -,line width=6pt}}]
				\node[shape=circle ,minimum size=0.2cm , draw=black,label={[blue]right:{$a$}}] (a) at (0,0){};
				\node[shape=circle,minimum size=0.2cm , draw=black] (b) at (1,1){};
				\node[shape=circle ,minimum size=0.2cm , draw=black] (c) at (-1,1){};
				\node[shape=circle ,minimum size=0.2cm , draw=black,label={[blue]right:{$b$}}] (d) at (0,2){};
				
				\path[-](a) edge node[auto]{}(b);
				\path[-](c) edge node[auto]{}(d);
				\path[-](b) edge node[auto]{}(d);
				\path[-](a) edge node[auto]{}(c);
			\end{tikzpicture}
		\end{minipage}			
	\end{center}
	\caption{Figures for the proof of Theorem \ref{EnumueratePlusOne} and \ref{EnumueratePlusTwo}.}\label{FigureProofFDLP1AndP2}
\end{figure}

In this section, we apply Theorem \ref{LiftingTheorem} to derive values of $\dede(n+1), \dede(n+2), \dede(n+3)$, by respectively counting all monotonous maps from $\mathbbm{2}^1,\mathbbm{2}^2$ and $\mathbbm{2}^3$ to $\mathbb{D}_n$. This is a warmup for the next section and provides an overview over formulas to compute some Dedekind numbers. 
\begin{theorem}\label{EnumueratePlusOne}
The following formulas compute $\mathbb{D}_{n+1}$ by means of $\mathbb{D}_{n}$:
\[\dede(n+1)=\sum_{a\in \dedebig(n)} \top(a)=\sum_{I\in\Int(\mathbb{D}_n)} 1 = \#\Int(\mathbb{D}_n) = \sum_{[I]\in\Int(\mathbbm{D}_n)/_\equiv} \#[I].\]
\begin{proof}
\emph{
If $\mathbbm{2}$'s bottom is mapped to $a\in\dedebig(n)$, every $b\in [a,\top]$ is a valid image for the top element. This implies that every  $a,b\in\dedebig(n)$, such that $[a,b]\in\Int(\mathbb{D}_n)$, is a valid choice too (see Figure \ref{FigureProofFDLP1AndP2}, left image). Instead of summing up all intervals, we can restrict the summation to equivalent ones. 
}\qed
\end{proof} 
\end{theorem}
The first equation from Theorem \ref{EnumueratePlusTwo} was presented in \cite{BermanKoehler76}.
\begin{theorem}\label{EnumueratePlusTwo}
	The following formulas enumerate $\mathbb{D}_{n+2}$ by means of $\mathbb{D}_{n}$:
	\[\dede(n+2)=\sum_{a,b\in \dedebig(n)} \bot(a\wedge b)\cdot\top(a\vee b)=\sum_{I\in\Int(\mathbb{D}_n)} (\#I)^2 = \sum_{[I]\in\Int(\mathbbm{D}_n)/_\equiv} (\#I)^2\cdot\#[I].\]
\begin{proof}
\emph{
If $a,b\in\dedebig(n)$ are chosen as in Figure \ref{FigureProofFDLP1AndP2}'s center, every element in $[a\vee b,\top]$ and each in $[\bot,a\wedge b]$
are possibilities for top and bottom respectively. On the other hand, Figure \ref{FigureProofFDLP1AndP2}'s right image shows that bottom and top can be mapped to every interval $[a,b]$, which implies the second last equation. Since isomorphic and anti isomorphic intervals have equal cardinality, the last equality holds.
}\qed
\end{proof}
\end{theorem}

The next theorem shows how to express the computation of $\dede(n+2)$ with matrix calculus. For that we use the trace operator $\Tr$ of square matrices, that is defined as the sum of all diagonal elements. 
\begin{theorem}\label{EnumueratePlusTwoMatrix}
	We define two matrices $\alpha(a,b):=\bot(a\wedge b)$ and $\beta(a,b):=\top(a\vee b)$, with $a,b\in \dedebig(n)$.
	It holds that: \[\dede(n+2)=\Tr(\alpha\cdot\beta).\]
	\begin{proof}
		\emph{
			Let $\gamma=\alpha\cdot\beta$ be the matrix product. We use the first formula from Theorem \ref{EnumueratePlusTwo} and insert the matrix expressions. Note that $\alpha$ and $\beta$ are symmetric.
			\[\dede(n+2)=\sum_{a\in \dedebig(n)}\sum_{b\in \dedebig(n)}\alpha(a,b)\cdot\beta(b,a)=\sum_{a\in \dedebig(n)}\gamma(a,a)=\Tr(\gamma).\]
		}\qed
	\end{proof}
\end{theorem}

\begin{figure}		
	\begin{center}
		\begin{minipage}{0.477\textwidth} 
		\begin{center}
			\begin{tikzpicture}[shorten >=0pt,shorten <=0pt,auto,cross line/.style={preaction={draw=white, -,line width=6pt}}]
				
				\node[shape=circle,minimum size=0.2cm , draw=black,label={[blue]right:{$a\wedge b\wedge c = x$}}] (1) at (0,0){};
				
				\node[shape=circle,minimum size=0.2cm , draw=black,label={[blue]right:{$c$}}] (2) at (1.5,1.5){};
				\node[shape=circle,minimum size=0.2cm , draw=black,label={[blue]right:{$b$}}] (3) at (0,1.5){};
				\node[shape=circle,minimum size=0.2cm , draw=black,label={[blue]right:{$a$}}] (4) at (-1.5,1.5){};
				
				\node[shape=circle,minimum size=0.2cm , draw=black,label={[blue]right:{$b\vee c$}}] (5) at (1.5,3){};
				\node[shape=circle,minimum size=0.2cm , draw=black,label={[blue]right:{$a\vee c$}}] (6) at (0,3){};
				\node[shape=circle,minimum size=0.2cm , draw=black,label={[blue]right:{$a\vee b$}}] (7) at (-1.5,3){};
				
				\node[shape=circle,minimum size=0.2cm , draw=black,label={[blue]right:{$y$}}] (8) at (0,4.5){};					
				
				\path[-](1) edge node[auto]{}(2);
				\path[-](1) edge node[auto]{}(3);
				\path[-](1) edge node[auto]{}(4);
				
				\path[-](8) edge node[auto]{}(5);
				\path[-](8) edge node[auto]{}(6);
				\path[-](8) edge node[auto]{}(7);
				
				\path[-](2) edge node[auto]{}(5);
				\path[-](2) edge node[auto]{}(6);
				
				\path[-](3) edge node[auto]{}(5);
				\path[-](3) edge node[auto]{}(7);
				
				\path[-](4) edge node[auto]{}(6);
				\path[-](4) edge node[auto]{}(7);
			\end{tikzpicture}
		\end{center}
		\end{minipage}
		\begin{minipage}{0.477\textwidth} 
		\begin{center}
			\begin{tikzpicture}[shorten >=0pt,shorten <=0pt,auto,cross line/.style={preaction={draw=white, -,line width=6pt}}]
				
				\node[shape=circle,minimum size=0.2cm , draw=black, label={[blue]right:{$x$}}] (1) at (0,0){};
				
				\node[shape=circle,minimum size=0.2cm , draw=black] (2) at (1.5,1.5){};
				\node[shape=circle,minimum size=0.2cm , draw=black, label={[blue]right:{$a$}}] (3) at (0,1.5){};
				\node[shape=circle,minimum size=0.2cm , draw=black, label={[blue]right:{$b\wedge c$}}] (4) at (-1.5,1.5){};
				
				\node[shape=circle,minimum size=0.2cm , draw=black] (5) at (1.5,3){};
				\node[shape=circle,minimum size=0.2cm , draw=black, label={[blue]right:{$b$}}] (6) at (0,3){};
				\node[shape=circle,minimum size=0.2cm , draw=black, label={[blue]right:{$c$}}] (7) at (-1.5,3){};
				
				\node[shape=circle,minimum size=0.2cm , draw=black, label={[blue]right:{$y$}}] (8) at (0,4.5){};					
				
				\path[dashed](1) edge node[auto]{}(2);
				\path[-](1) edge node[auto]{}(3);
				\path[-](1) edge node[auto]{}(4);
				
				\path[dashed](8) edge node[auto]{}(5);
				\path[-](8) edge node[auto]{}(6);
				\path[-](8) edge node[auto]{}(7);
				
				\path[dashed](2) edge node[auto]{}(5);
				\path[dashed](2) edge node[auto]{}(6);
				
				\path[dashed](3) edge node[auto]{}(5);
				\path[-](3) edge node[auto]{}(7);
				
				\path[-](4) edge node[auto]{}(6);
				\path[-](4) edge node[auto]{}(7);
			\end{tikzpicture}
		\end{center}
		\end{minipage}			
	\end{center}
	\caption{Figures for the proof of Theorem \ref{EnumueratePlusThree}.}\label{ProofFDLPlusThree}
\end{figure}
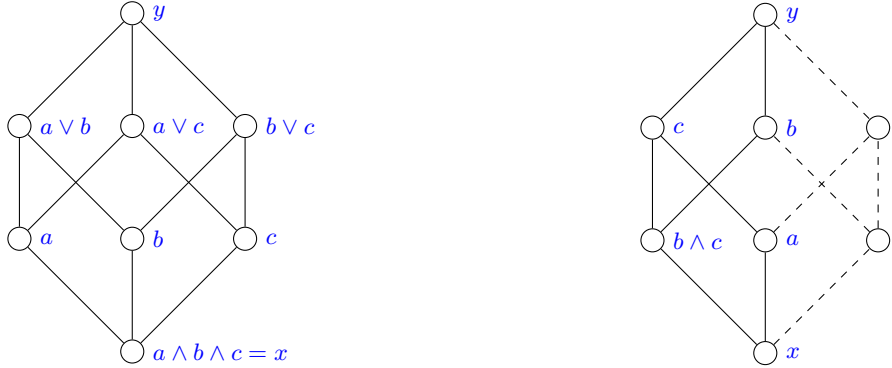

The next theorem's first formula is similar to the one stated in \cite{Hoedt15}.
\begin{theorem}\label{EnumueratePlusThree}
Let $X$ denote the set of all intervals $\{[\bot,y]\mid y\in \dedebig(n)\}$.
The following formulas enumerate $\mathbb{D}_{n+3}$ by means of $\mathbb{D}_{n}$:
\begin{align*}
\dede(n+3)&=\sum_{y\in \dedebig(n)}\sum_{a,b,c\in [\bot,y]}\bot(a\wedge b\wedge c)\cdot\top(a\vee b)\cdot\top(a\vee c)\cdot\top(b\vee c)\\
		  &=\sum_{[I]\in X/_{\cong}} \#[I]\cdot\sum_{a,b,c\in I} \bot(a\wedge b\wedge c)\cdot\top(a\vee b)\cdot\top(a\vee c)\cdot\top(b\vee c)\\
		  &=\sum_{I\in\Int(\mathbb{D}_n)} \sum_{a,b,c\in I}\top(a\vee b)\cdot\top(a\vee c)\cdot\top(b\vee c)\\
		  &=\sum_{[I]\in\Int(\mathbbm{D}_n)/_\equiv} \#[I]\cdot\sum_{a,b,c\in I}\top(a\vee b)\cdot\top(a\vee c)\cdot\top(b\vee c)\\
		  &=\sum_{[I]\in\Int(\mathbbm{D}_n)/_\equiv} \#[I]\cdot\sum_{a,b\in I}(\sum_{c\in [a,\top_{\I}]}\bot(b\wedge c))^2.
\end{align*}
\begin{proof}
\emph{
Choosing $y$ and then $a,b,c\in [\bot,y]$ as in Figure \ref{ProofFDLPlusThree}'s left image, the first formula can be deduced. By factoring out isomorphic elements from $X$, we derive the second equation. It is important to notice that after choosing $a,b,c\in [\bot,y]$, there is still one degree of freedom for the bottom element. If we map bottom and top to an interval $[x,y]$'s boundaries, we remove this degree of freedom and get the third equation. Obviously, isomorphic intervals can then be factored out. Equation number four follows from the fact that anti isomorphic intervals can be factored out as well. To see why this holds, we dualize the third identity\footnote{This implies $a,b,c$ have to be placed one "level" above than depicted in Figure \ref{ProofFDLPlusThree}.}:
\[\dede(n+3) = \sum_{[I]\in\Int(\mathbb{D}_n)} \sum_{a,b,c\in I}\bot(a\wedge b)\cdot\bot(a\wedge c)\cdot\bot(b\wedge c).\] 
Since an anti isomorphism reverses order and thereby join and meet operations, for $I\conga J$ it holds that:
\[\sum_{a,b,c\in I}\bot(a\wedge b)\cdot\bot(a\wedge c)\cdot\bot(b\wedge c) = \sum_{a,b,c\in J}\top(a\vee b)\cdot\top(a\vee c)\cdot\top(b\vee c).\]
To conclude the last equation, we utilize that $\mathbbm{2}^3\cong \mathbbm{2}^2 \times\mathbbm{2}$. The goal is to apply vertical symmetry from $\mathbbm{2}^3$ (see Figure \ref{ProofFDLPlusThree}, right image). If an interval $I$ is chosen, every $a,b\in I=[x,y]$ can be independently placed as depicted. After computing $\sum_{c\in [a,y]}\#[x,b\wedge c]$, the result can be squared due to symmetry induced by $\mathbbm{2}^2$. 
}\qed
\end{proof}
\end{theorem}

Like in Theorem \ref{EnumueratePlusTwoMatrix}, we express the computation of $\dede(n+3)$ via matrix calculus.
\begin{theorem}\label{EnumueratePlusThreeMatrix}
	We define the matrix $\alpha(a,b):=\top(a\vee b)$, with $a,b\in I\in\Int(\mathbbm{D}_n)$.
	It holds that: \[\dede(n+3)=\sum_{[I]\in \Int(\mathbbm{D}_n)/_\equiv} \#[I]\cdot\Tr(\alpha^3).\]
	\begin{proof}
		\emph{
			We use Theorem \ref{EnumueratePlusTwo}'s fourth formula and insert the matrix definition.
			\[\dede(n+3)=\sum_{[I]\in \Int(\mathbbm{D}_n)/_\equiv} \#[I]\cdot \sum_{a\in I}\sum_{b\in I}\alpha(a,b)\sum_{c\in I}\alpha(b,c)\alpha(c,a).\]
			The innermost sum expresses $\alpha^2$, the sum over $b$ leads to $\alpha^3$ and the sum over $a$ gives the trace operator. 
		}\qed
	\end{proof}
\end{theorem}

Our implementation of Theorem \ref{EnumueratePlusThreeMatrix} with GPU support, run on a Nvidia A100 GPU, computes $\dede(8)$ in $17.56$ seconds. However, this algorithm can not compute $\dede(9)$, since the involved matrices would have a maximal dimension of $7828354\times 7828354$.

\section{Enumeration of $\mathbb{D}_{n+4}$}\label{P4Section}
In this section, we apply Theorem \ref{LiftingTheorem} to derive values of $\dede(n+4)$, by counting all monotonous maps from $\mathbbm{2}^4$ to $\mathbb{D}_n$.
\begin{theorem}\label{EnumueratePlusFour}	
	The following formulas enumerate $\mathbb{D}_{n+4}$ by means of $\mathbb{D}_{n}$:
	\[\dede(n+4)=\sum_{I\in\Int(\mathbb{D}_n)} X = \sum_{[I]\in\Int(\mathbbm{D}_n)/_\equiv} \#[I]\cdot X,\;\text{with:}\]
	\begin{align*}
		X&=\sum_{a,b,c,d\in I}\sum_{\substack{e\in[a\vee b,\top] \\ f\in[a\vee c,\top]\\ g\in[a\vee d,\top]\\ h\in[b\vee c,\top]\\ i\in[b\vee d,\top]\\ j\in[c\vee d,\top]}}\top(e\vee f\vee h)\cdot\top(e\vee g\vee i)\cdot\top(f\vee g\vee j)\cdot\top(h\vee i\vee j)\\
		&=\sum_{a,b,c,d\in I}\sum_{\substack{e\in[a\vee b\vee c,\top] \\ f\in[a\vee b\vee d,\top]\\ g\in[a\vee c\vee d,\top]\\ h\in[b\vee c\vee d,\top]}}(\;\#[a\vee b,e\wedge f]\cdot\#[a \vee c,e\wedge g]\cdot\#[a \vee d,f\wedge g]\\
		&\;\;\;\;\;\;\;\;\;\;\;\;\;\;\;\;\;\;\;\;\;\;\;\;\;\;\;\;\;\;\;\;\;\cdot\#[b \vee c,e\wedge h]\cdot\#[b\vee d,f\wedge h]\cdot\#[c\vee d,g\wedge h]\;)\\
		&=\sum_{a,b,c,d,e,f\in I}(\bot(a\wedge c\wedge d)\cdot\bot(b\wedge c\wedge e)\cdot\bot(b\wedge d\wedge f)\cdot\bot(a\wedge e\wedge f)\\ &\;\;\;\;\;\;\;\;\;\;\;\;\;\;\;\;\;\;\;\;\cdot\top(b\vee c\vee d)\cdot\top(a\vee c\vee e)\cdot\top(a\vee d\vee f)\cdot\top(b\vee e\vee f)).
	\end{align*}
	\begin{proof}
		\emph{
		For all three formulas, bottom and top of $\mathbbm{2}^4$ are mapped to the bounds of every possible interval $[x,y]$. But as in Theorem \ref{EnumueratePlusThree}'s proof, we can argue that it is sufficient to consider equivalent intervals.\\		
		In the first equation, values $a,b,c,d$ are chosen as depicted by the red labels of Figure \ref{TwoTimesTwoFigure}. We then iterate over all six ranges given by the $2$-joins and the interval top. This leaves four degrees of freedom which range form the respective $3$-joins up to the top element.\\
		The second equation uses the same choice for $a,b,c,d$. In the following step, we iterate over all four ranges given by the $3$-joins and the top element via $e,f,g,h$. This results in six degrees of freedom that range from the $2$-joins of $a,b,c,d$ up to the respective $2$-meet of $e,f,g,h$.\\
		Lastly, we can place $a,b,c,d,e,f$ as depicted by the blue labels of Figure \ref{TwoTimesTwoFigure}. Consequently, there are eight degrees of freedom that range from every $3$-join to the interval top, and from the interval bottom to every $3$-meet.
		}\qed	
\begin{figure}[h]
			\begin{center}
				\scalebox{.8}{
				\begin{tikzpicture}[scale=0.85,shorten >=0pt,shorten <=0pt,auto,cross line/.style={preaction={draw=white, -,line width=6pt}}]
					
					\node[shape=circle,minimum size=0.2cm , draw=black, label={[blue]below:{$x$}}] (1) at (0,-0.5){};					
					\node[shape=circle,minimum size=0.2cm , draw=black, label={[red]right:{$b$}}, label={[blue]left:{$b\wedge c\wedge e$}}] (2) at (-1,2){};
					\node[shape=circle,minimum size=0.2cm , draw=black, label={[red]left:{$c$}}, label={[blue]right:{$b\wedge d\wedge f$}}] (3) at (1,2){};				
					\node[shape=circle,minimum size=0.2cm , draw=black, label={[blue]above:{$b$}}] (4) at (0,4.5){};					
					
					\node[shape=circle,minimum size=0.2cm , draw=black,label={[red]right:{$a$}}, label={[blue]left:{$a\wedge c\wedge d$}}] (5) at (-5,3.5){};					
					\node[shape=circle,minimum size=0.2cm , draw=black, label={[blue]left:{$c$}}] (6) at (-6,6){};
					\node[shape=circle,minimum size=0.2cm , draw=black, label={[blue]right:{$d$}}] (7) at (-4,6){};
					\node[shape=circle,minimum size=0.2cm , draw=black, label={[blue]left:{$b\vee c\vee d$}}] (8) at (-5,8.5){};
					
					\node[shape=circle,minimum size=0.2cm , draw=black, label={[red]left:{$d$}}, label={[blue]right:{$a\wedge e\wedge f$}}] (9) at (5,3.5){};
					\node[shape=circle,minimum size=0.2cm , draw=black, label={[blue]right:{$f$}}] (10) at (6,6){};
					\node[shape=circle,minimum size=0.2cm , draw=black, label={[blue]left:{$e$}}] (11) at (4,6){};					
					\node[shape=circle,minimum size=0.2cm , draw=black, label={[blue]right:{$b\vee e\vee f$}}] (12) at (5,8.5){};
					
					\node[shape=circle,minimum size=0.2cm , draw=black, label={[blue]below:{$a$}}] (13) at (0,7.5){};				
					\node[shape=circle,minimum size=0.2cm , draw=black, label={[blue]left:{$a\vee c\vee e$}}] (14) at (-1,10){};
					\node[shape=circle,minimum size=0.2cm , draw=black, label={[blue]right:{$a\vee d\vee f$}}] (15) at (1,10){};					
					\node[shape=circle,minimum size=0.2cm , draw=black, label={[blue]above:{$y$}}] (16) at (0,12.5){};	
					
					\path[-](1) edge node[auto]{}(2);
					\path[-](1) edge node[auto]{}(3);
					\path[-](2) edge node[auto]{}(4);
					\path[-](3) edge node[auto]{}(4);
					
					\path[-](16) edge node[auto]{}(15);
					\path[-](16) edge node[auto]{}(14);
					\path[-](15) edge node[auto]{}(13);
					\path[-](14) edge node[auto]{}(13);
					
					\path[-](5) edge node[auto]{}(6);
					\path[-](5) edge node[auto]{}(7);
					\path[-](8) edge node[auto]{}(6);
					\path[-](8) edge node[auto]{}(7);
					
					\path[-](12) edge node[auto]{}(11);
					\path[-](12) edge node[auto]{}(10);
					\path[-](9) edge node[auto]{}(11);
					\path[-](9) edge node[auto]{}(10);
					
					\path[-](1) edge node[auto]{}(5);
					\path[-](5) edge node[auto]{}(13);
					
					\path[-](4) edge node[auto]{}(8);
					\path[-](8) edge node[auto]{}(16);
					
					\path[-](2) edge node[auto]{}(6);
					\path[-](6) edge node[auto]{}(14);
					
					\path[-](3) edge node[auto]{}(7);
					\path[-](7) edge node[auto]{}(15);
					
					\path[-](1) edge node[auto]{}(9);
					\path[-](9) edge node[auto]{}(13);
					
					\path[-](4) edge node[auto]{}(12);
					\path[-](12) edge node[auto]{}(16);
					
					\path[-](2) edge node[auto]{}(11);
					\path[-](11) edge node[auto]{}(14);
					
					\path[-](3) edge node[auto]{}(10);
					\path[-](10) edge node[auto]{}(15);					
				\end{tikzpicture}
			}
	\end{center}
	\caption{Figure of $\mathbbm{2}^4$ seen as $\mathbbm{2}^2 \times \mathbbm{2}^2$.}\label{TwoTimesTwoFigure}
\end{figure}
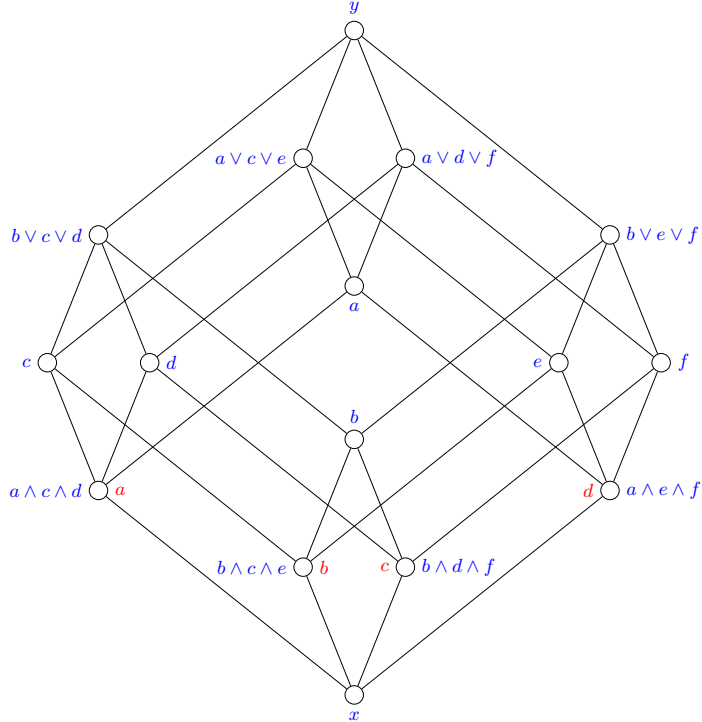
\end{proof}
\end{theorem}

Theorem \ref{EnumueratePlusFour}'s last formula can be expressed via matrix calculus.
\begin{theorem}\label{MatrixTheorem}
For $I\in\Int(\mathbbm{D}_n)$ and $a,b,c,d\in I$, we define matrices $\alpha$ and $\beta$:
\[\alpha_{ab}(c,d):=\bot(a\wedge c\wedge d)\cdot\top(b\vee c\vee d)\;\text{and}\;\beta_{ab}(c,d):=\bot(b\wedge c\wedge d)\cdot\top(a\vee c\vee d).\]
Let $\gamma_{ab}$ be the matrix product of $\alpha_{ab}$ and $\beta_{ab}$. It holds that:
\[\dede(n+4)=\sum_{[I]\in\Int(\mathbbm{D}_n)/_\equiv} \#[I]\cdot \sum_{a,b\in I} \Tr(\gamma^2_{ab}).\]
\begin{proof}
\emph{
	Inserting the matrix expressions into the formula from Theorem \ref{EnumueratePlusFour}, we get for $d(n+4):$
	\[\sum_{[I]\in\Int(\mathbbm{D}_n)/_\equiv} \#[I]\cdot \sum_{a,b\in I}\sum_{c,d\in I}\sum_{e,f\in I}\alpha(c,d)\beta(c,e)\alpha(e,f)\beta(d,f).\]
	\phantom{Mi}Since $\alpha$ and $\beta$ are symmetric, we can identify two matrix multiplications: \[\gamma(d,e)=\sum_{c\in I}\alpha(d,c)\beta(c,e)\;\text{and}\;\gamma(e,d)=\sum_{f\in I}\alpha(e,f)\beta(f,d).\] 
	\phantom{Mi}All indices are from the same range (the interval $I$), which means that the matrices above are equal. Proceeding from this matrix point of view, we compute the diagonal of the product from $\gamma$ with itself and sum it up $\sum_{d\in I}\sum_{e\in I}\gamma(d,e)\gamma(e,d)$.}\qed
\end{proof}
\end{theorem}
Looking at Figure \ref{TwoTimesTwoFigure}, we see that $a$ and $b$ are independent from each other. Together with the formula from Theorem \ref{EnumueratePlusFour}, this implies rotation symmetry between $a$ and $b$. It is enough to consider pairs such that $a\leq b$. Values for $a$ that are strictly smaller than $b$ are weighted with factor two and for $a=b$ the weight is one. This reduces computations by almost one half. Still, further reductions are possible. Therefor, for every $I\in\Int(\mathbbm{D}_n)$, let $\varphi: I\rightarrow I$ be an (anti) isomorphism. We introduce a relation on $I\times I$. That is $(a,b)\sim (\tilde{a},\tilde{b}):\Longleftrightarrow$
\[\exists \varphi: I\rightarrow I,\;\varphi(a)=\varphi(\tilde{a})\;\text{and}\;\varphi(b)=\varphi(\tilde{b}).\]

Let $I\times I\mid_{\leq}$ denote all pairs $(a,b)$, such that $a\leq b$. Furthermore the weight operator $\omega: I\times I\mid_{\leq}\,\rightarrow \{1,2\}$ equals $2$ for $a\neq b$ and $1$ otherwise.
\begin{theorem}\label{MatrixTheorem2}
	For $I\in\Int(\mathbbm{D}_n)$ and $a,b,c,d\in I$, we define matrices $\alpha$ and $\beta$:
	\[\alpha_{ab}(c,d):=\bot(a\wedge c\wedge d)\cdot\top(b\vee c\vee d)\;\text{and}\;\beta_{ab}(c,d):=\bot(b\wedge c\wedge d)\cdot\top(a\vee c\vee d).\]
	Let $\gamma_{ab}$ be the matrix product of $\alpha_{ab}$ and $\beta_{ab}$. It holds that:
	\[\dede(n+4)=\sum_{[I]\in\Int(\mathbbm{D}_n)/_\equiv} \#[I]\cdot \sum_{[(a,b)]\in (I\times I\mid_{\leq})/_\sim} \omega(a,b)\cdot \#[(a,b)]\cdot\Tr(\gamma^2_{ab}).\]
	\begin{proof}
		\emph{
			The relation $\sim$ is an equivalence relation and all pairs within one equivalence class result in an equal summand for the enumeration formula. We illustrate this at the simpler example $\sum_{z\in I} \bot(a\wedge z)\cdot\top(a\vee z)$. Under an isomorphism $\varphi$, with $\varphi(a)=\tilde{a}$ (also note that top and bottom are preserved), we get: \[\sum_{\varphi(z)\in I} \bot(\varphi(a\wedge z))\cdot\top(\varphi(a\vee z))=\sum_{\varphi(z)\in I} \bot(\tilde{a}\wedge \varphi(z))\cdot\top(\tilde{a}\vee \varphi(z))=\sum_{\tilde{z}\in I} \bot(\tilde{a}\wedge \tilde{z})\cdot\top(\tilde{a}\vee \tilde{z}).\] Since $\varphi$ is bijective and preserves join, meet and order, both sums (with and without applying $\varphi$) are equal. An anti isomorphism would inverse order, but this is offset by the symmetry w.r.t. top and bottom. These arguments can be extended to the actual formula from Theorem \ref{EnumueratePlusFour}.\qed			
		}
	\end{proof}
\end{theorem}
Our implementation of Theorem \ref{MatrixTheorem2} with GPU support, as described in the following section, can compute $\dede(8)$ in about $3$ seconds.
\section{Algorithm To Compute The Ninth Dedekind Number}\label{NinthDedekindAlgo}
We list all steps necessary to compute $\dede(n+4)$ via Theorem \ref{MatrixTheorem2} and address them individually. 
\begin{enumerate}
	\item Generate $\mathbbm{D}_n$ as described at the end of Section \ref{Intro}.
	\item Compute the equivalence classes $\Int(\mathbbm{D}_n)/_\equiv$ and save one representative for each class together with the class's cardinality.
	\item For every equivalence class representative $I$ of $\Int(\mathbbm{D}_n)/_\equiv$, compute equivalence classes of $(I\times I\mid_{\leq})/_\sim$. Again save one representative for each class together with the respective cardinality. 
	\item Run the computation according to Theorem \ref{MatrixTheorem2}:
	\begin{algorithm*}	
		\KwData{$\Int(\mathbbm{D_n})/_\equiv$ and $\forall[I]\,\in\Int(\mathbbm{D_n})/_\equiv:\;(I\times I\mid_{\leq})/_\sim$}
		\KwResult{$\dede(n+4)$}	
		\For{$[I]\in\Int(\mathbbm{D_n})/_\equiv$}
		{		
			\For{$[(a,b)]\in(I\times I\mid_{\leq})/_\sim$}
			{
				- generate the matrices  $\alpha_{ab}$ and $\beta_{ab}$;\\
				- compute the matrix product $\gamma_{ab}=\alpha_{ab}\cdot\beta_{ab}$;\\
				- compute the trace of $\gamma^2_{ab}$;\\
				- multiply the trace with $\omega(a,b)$ and $\#[(a,b)]$;\\
			}
			- sum up each value from above;\\
			- multiply the sum with $\#[I]$;\\
		}	
	\end{algorithm*}
\end{enumerate}
\subsection{Step 1}
The generation of $\mathbbm{D}_n$ is straight forward and needs no special explanation. Elements of $\mathbbm{D}_4$ can be represented by $16$-bit, $\mathbbm{D}_5$ by $32$-bit and $\mathbbm{D}_6$ by $64$-bit unsigned integers. To store all elements of $\mathbbm{D}_6$ takes about $63$ MB. Even the elements of $\mathbbm{D}_7$ can be generated within a few hours, but to store them would require about $38.63$ TB of memory, which makes $\mathbbm{D}_7$ not suitable as a basis to calculate $\mathbbm{D}_9$.
\subsection{Step 2}
It is a simple task to generate all intervals $\Int(\mathbbm{D}_n)$. In the following, we describe how to determine equivalence classes w.r.t. $\equiv$, like depicted in Example \ref{IsoIntervalExample}. For that we use \emph{Formal Concept Analysis} (see \cite{GanterWille}), which treats the relationship between complete lattices and binary relations. We present some aspects of this theory that we will utilize to achieve our goals.\\

\phantom{MI}Given a complete lattice $\mathbb{L}=(L,\leq)$, a \emph{join irreducible} (\emph{meet irreducible}) element of $\mathbb{L}$ can not be expressed as a proper join (meet) of other elements form $\mathbb{L}$ that are strictly below (above). Let $J(\mathbb{L})$ denote all join irreducibles and $M(\mathbb{L})$ all meet irreducibles. It holds that $\mathbb{L}$ can be represented by the binary relation $\mathbb{R}:=(J(\mathbb{L}),M(\mathbb{L}),R)$ with	
$R:=\leq\mid_{J(\mathbb{L})\times M(\mathbb{L})}$, that is referred to as a \emph{formal context}. Every element, as well as various properties and parameters of $\mathbb{L}$ can be reconstructed from the formal context. For instance, every $x\in L$ can be represented by the set $J(x)$ (all join irreducibles equal to or below $x$) or $M(x)$ (all meet irreducibles equal to or above $x$). Since the transition form $\mathbb{L}$ to $\mathbb{R}$ constitutes an up to logarithmic size reduction, it is often beneficial to use $\mathbb{R}$ for computations.\\

\phantom{MI}It is shown in \cite{GanterWille} that an interval $[x,y]\in \Int(\mathbbm{L})$, can be described by a subrelation of $R$ namely $(J(y),M(x),R\cap(J(y)\times M(x)))$. Furthermore, a lattice isomorphism $\varphi: \mathbb{L}_1\rightarrow\mathbb{L}_2$ exists if and only if a \emph{context isomorphism} $(\alpha,\beta):\mathbb{R}_1\rightarrow \mathbb{R}_2$ exists. The latter one is defined as a pair of bijective maps $\alpha: J(\mathbb{R}_1)\rightarrow J(\mathbb{R}_2)$ and $\beta: M(\mathbb{R}_1)\rightarrow M(\mathbb{R}_2)$, such that $(j,m)\in R_1 \Longleftrightarrow (\alpha(j),\beta(m))\in R_2$. A lattice anti isomorphism corresponds to a context isomorphism from $\mathbb{R}_1$ to the \emph{dual context} $\mathbb{R}^d_2:=(M(\mathbb{R}_2),J(\mathbb{R}_2),R_2^d)$, with $R^d_2:=\geq\cap (M(\mathbb{R}_2)\times J(\mathbb{R}_2))$.\\

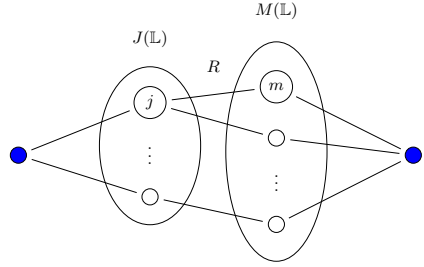
\begin{figure}[h]
	\begin{center}
		\scalebox{.7}{
			\begin{tikzpicture}[every node/.style={draw,circle},fsnode/.style={fill=white},ssnode/.style={fill=white},bbnode/.style={fill=blue},
				every fit/.style={ellipse,draw,inner sep=-2pt,text width=1.5cm},shorten >= 3pt,shorten <= 3pt,scale=0.6]				
				\begin{scope}[start chain=going below,node distance=5mm]
					\node[fsnode,on chain] (f1){$j$};
					\node[white,on chain,label={[xshift=0cm, yshift=-0.5cm]$\vdots$}] (f2){};
					\node[fsnode,on chain] (f3){};
				\end{scope}				
				\begin{scope}[xshift=4cm,yshift=0.5cm,start chain=going below,node distance=5mm]
					\node[ssnode,on chain] (s4){$m$};
					\node[ssnode,on chain] (s5){};
					\node[white,on chain,label={[xshift=0cm, yshift=-0.5cm]$\vdots$}] (s6){};
					\node[ssnode,on chain] (s7){};					
				\end{scope}				
				
				\node [black,fit=(f1) (f3),label=above:$J(\mathbb{L})$] {};				
				\node [black,fit=(s4) (s7),label=above:$M(\mathbb{L})$] {};							
				\node [white,yshift=0.2cm,xshift=1.2cm,label=above:{$R$}] {};
				\node[bbnode,yshift=-1cm,xshift=5cm] (s8){};
				\node[bbnode,yshift=-1cm,xshift=-2.5cm] (f4){};
				
				\draw (f4) -- (f1);
				\draw (f4) -- (f3);				
				\draw (s8) -- (s4);
				\draw (s8) -- (s5);
				\draw (s8) -- (s7);
				
				\draw (f1) -- (s5);
				\draw (f1) -- (s4);				
				\draw (s7) -- (f3);
			\end{tikzpicture}
		}
	\end{center}
	\caption{A formal context interpreted as bipartite graph.}\label{ContextAsGraphFigure}
\end{figure}

\phantom{MI}The identification of (anti) isomorphic formal contexts can be translated to an isomorphism problem of colored bipartite graphs. That is because every formal context represents a bipartite graph (see Figure \ref{ContextAsGraphFigure}). If we want to determine context isomorphisms only, each of the two vertex sets gets a different color label assigned to, and we compute graph isomorphisms that respect them. Since we want to determine anit isomorphisms too, we have to allow for a complete switch of the vertex sets, but still no mapping from one set into the other. This can be achieved by adding two additional equally colored vertices, and connecting them with every vertex of the bipartite set that they represent.

To practically determine isomorphic graphs, we use the software "nauty" (see \cite{NautyTraces}), that can compute a canonical string representation of a given colored graph. This translates the graph isomorphism problem into a string comparison, since two graphs are isomorphic if and only if their canonical string representations are equal.\\


%
%

	\begin{algorithm*}
		\KwData{$\Int(\mathbbm{D_n})$}
		\KwResult{$\Int(\mathbbm{D_n}/_\equiv$}	
		\For{$[x,y]\in\Int(\mathbbm{D_n})$}
		{
			- compute a formal context that represents $[x,y]$;\\
			- transform the context to a colored bipartite graph as depicted in Figure \ref{ContextAsGraphFigure};\\ %
			- use "nauty" to compute a canonical string;\\
			- count occurrences of each string;\\
		}	
		\caption{The steps to calculate $\Int(\mathbbm{D_n})/_\equiv$.}\label{IntervalAlgorithm}
	\end{algorithm*}

For $\mathbbm{D}_2$ we compute 6, for $\mathbbm{D}_3$ 18, for $\mathbbm{D}_4$ 134 and for $\mathbbm{D}_5$ 9919 equivalence classes. Due to nauty's efficiency, the last computation can be performed within just a couple of seconds on a $2.4$ GHz desktop PC. Note that according to Theorem \ref{EnumueratePlusOne} there are $7828354$ intervals of $\mathbbm{D}_5$. Hence we reduce them to just $0.00127 \%$, which has a huge runtime impact on the computation of $\dede(5+4)$. 
%
\subsection{Step 3}
\begin{figure}
	\begin{center}
		\scalebox{.7}{
			\begin{tikzpicture}[every node/.style={draw,circle},fsnode/.style={fill=white},ssnode/.style={fill=white},bbnode/.style={fill=blue},rrnode/.style={fill=red},
				every fit/.style={ellipse,draw,inner sep=-2pt,text width=1.5cm},shorten >= 3pt,shorten <= 3pt,scale=0.6]				
				\begin{scope}[start chain=going below,node distance=5mm]
					\node[fsnode,on chain] (f1){$j$};
					\node[white,on chain,label={[xshift=0cm, yshift=-0.5cm]$\vdots$}] (f2){};
					\node[fsnode,on chain] (f3){};
				\end{scope}				
				\begin{scope}[xshift=4cm,yshift=0.5cm,start chain=going below,node distance=5mm]
					\node[ssnode,on chain] (s4){$m$};
					\node[ssnode,on chain] (s5){};
					\node[white,on chain,label={[xshift=0cm, yshift=-0.5cm]$\vdots$}] (s6){};
					\node[ssnode,on chain] (s7){};					
				\end{scope}				
				
				\node [black,fit=(f1) (f3),label=above:$J(\mathbb{L})$] {};				
				\node [black,fit=(s4) (s7),label=above:$M(\mathbb{L})$] {};							
				\node [white,yshift=0.2cm,xshift=1.2cm,label=above:{$R$}] {};
				\node[bbnode,yshift=-1cm,xshift=5cm] (s8){};
				\node[bbnode,yshift=-1cm,xshift=-2.5cm] (f4){};
				
				\node[fsnode,yshift=1cm,xshift=5cm,label=right:$m(a)$] (s9){};
				\node[fsnode,yshift=1cm,xshift=-2.5cm,label=left:$j(a)$] (f5){};
				
				\node[fsnode,yshift=-3cm,xshift=5cm,label=right:$m(b)$] (s10){};
				\node[fsnode,yshift=-3cm,xshift=-2.5cm,label=left:$j(b)$] (f6){};
				
				\node [rrnode,yshift=3cm,xshift=1.2cm] (x){};
				\node [rrnode,yshift=-3.5cm,xshift=1.2cm](y){};
				
				\draw[shorten >=0.5cm] (f5) -- (f1);
				\draw[shorten >=1.1cm] (f5) -- (f2);
				\draw[shorten >=1cm] (f6) -- (f2);
				\draw[shorten >=0.55cm] (f6) -- (f3);
				
				\draw[shorten >=0.5cm] (s9) -- (s4);
				\draw[shorten >=1cm] (s9) -- (s5);
				\draw[shorten >=0.9cm] (s10) -- (s6);
				\draw[shorten >=0.6cm] (s10) -- (s7);
				
				\draw (x) -- (s9);
				\draw (x) -- (f5);
				\draw (y) -- (s10);
				\draw (y) -- (f6);
				
				\draw (f4) -- (f1);
				\draw (f4) -- (f3);				
				\draw (s8) -- (s4);
				\draw (s8) -- (s5);
				\draw (s8) -- (s7);
				
				\draw (f1) -- (s5);
				\draw (f1) -- (s4);				
				\draw (s7) -- (f3);
			\end{tikzpicture}
		}
	\end{center}
	\caption{The graph that decodes the (anti) isomorphism problem w.r.t. $\sim$.}\label{IsoConnfigurationFigure}
\end{figure}
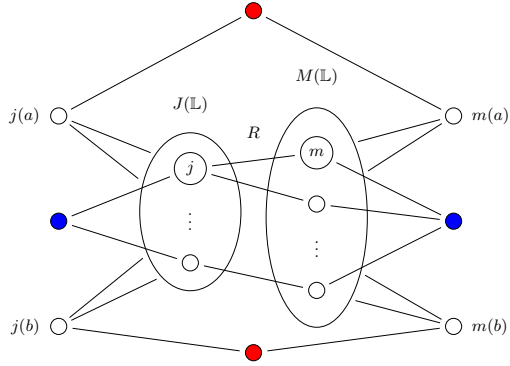

Considering an interval $I\in\Int(\mathbbm{D_n})$, we can generate all pairs $(a,b)\in I\times I\mid_\leq$. In order to determine the existence of an (anti) isomorphism $\varphi: I\rightarrow I$, with $\varphi(a)=\tilde{a}$ and $\varphi(b)=\tilde{b}$, for given pairs $(a,b)$ and $(\tilde{a},\tilde{b})$, we use a refined version of the approach from Step 2, which is illustrated in Figure \ref{IsoConnfigurationFigure}.\\

Firstly, the bipartite graph that decodes the (anti) isomorphism problem on $I$ (as in Step 2) has to be generated. Next we impose further restrictions, by adjoining 4 additional vertices that are connected to all join and meet irreducibles w.r.t. $a$ and $b$ respectively. The red vertices in Figure \ref{IsoConnfigurationFigure} connect the vertices related to $a$ and $b$ separately. This configuration assures that, if we compute equivalence classes w.r.t. nauty's canonical string representation, an (anti) isomorphism on $I$ respects the pair $(a,b)$, and can additionally swap $a$ and $b$. Example \ref{IsoConfigurationExample} shows the equivalence classes of $([0,15]\times[0,15])\mid_\leq$.

\begin{example}\label{IsoConfigurationExample}
	Equivalence classes of $([0,15]\times[0,15])\mid_\leq$ w.r.t. $\sim$.
	\begin{align*}
		& \{(0,0),(15,15)\}, \{(1,1),(7,7)\}, \{(3,3),(5,5)\}, \{(0,1),(7,15)\}, \{(3,5)\}, \{(0,3),\\
		& (0,5),(3,15),(5,15)\},\{(1,3),(1,5),(3,7),(5,7)\},\{(1,7)\},\{(0,7),(1,15)\},\{(0,15)\}.
	\end{align*}  
\end{example}

If we perform the computation for every $[I]\in \Int(\mathbbm{D_n})/_\equiv$, we treat a certain amount of pairs $(a,b)$ and reduce them to equivalence classes. We get the reductions $\mathbbm{D}_2:\,56 \rightarrow 33$, $\mathbbm{D}_3:\,1127 \rightarrow 446$, $\mathbbm{D}_4:\,274409 \rightarrow 80741$ and $\mathbbm{D}_5:\,8646896880 \rightarrow 4257682565$. The computation for $\mathbbm{D}_5$ took about one hour on an AWS "c6i.32xlarge" cluster with 128 cores. One could think that it is not worth the effort, since the reduction for $\mathbbm{D}_5$ is only slightly more than a factor two. But it can be observed that large intervals usually have a higher reduction rate. That is because they leave "more room" for symmetry than small intervals. For instance the largest interval of $\mathbbm{D}_5$ is $[0,4294967295]$. For this interval we have a reduction $57471561\rightarrow 140736$. This reduction has a big runtime impact, since the involved matrices have dimension $7581\times 7581$.
\subsection{Step 4}
For a given $I\in \Int(\mathbbm{D_n})$, we precompute all values of $\bot(\cdot)$ and $\top(\cdot)$ on a CPU host and transfer them to a GPU device. There the matrices $\alpha$ and $\beta$ are generated and multiplied. A batch of matrices is processed in parallel via CUDA's "cublasDgemmStridedBatched" kernel. Since a shift from integer values to double precision occurs, we have to assure that all values are not bigger than $2^{53}$ to rule out a loss of precision. Actually, we use the smaller bound $2^{51}$, since this is beneficial for another estimation.\\

An upper bound for the entries of $\alpha$ and $\beta$ is given by $\max_\alpha=\bot(a)\cdot\top(b)$ and $\max_\beta=\bot(b)\cdot\top(a)$. This leads to  $\max_\gamma=\max_\alpha\cdot\max_\beta\cdot\#I$. We compute $\max_\gamma$ for all pairs $(a,b)$ w.r.t. $\mathbbm{D}_5$. In $118084$ cases this estimation exceeds $2^{51}$. For these we perform the exact matrix multiplication with $64$-bit integer precision using the C++ Eigen library, and confirm that $\max_\gamma$ is smaller than $2^{51}$.\\

After multiplying the matrices, a handwritten CUDA kernel computes the trace of $\gamma^2$. For that $128$-bit unsigned integer values are used. Since the largest entry in $\gamma$ is smaller than $2^{51}$ and $7581\leq 2^{13}$, we can estimate a maximal value for the trace of $\gamma^2$ via $(2^{51} \cdot 2^{51} \cdot 2^{13}) \cdot 2^{13} = 2^{128}$. This assures us that $128$-bit unsigned integer precision is enough to compute the trace. Finally, the summation of all trace values and multiplication with the equivalence classes' carnality is done via $1024$-bit unsigned integers on the host system.  

\section{Conclusion}
We gave an overview about formulas to compute Dedekind numbers in Section \ref{WarmupSection} and \ref{P4Section}. Some of these formulas can by interpreted in terms of matrix multiplication, which allows for fast "number crunching" on a GPU device. In Section \ref{NinthDedekindAlgo}, the steps of an algorithm to compute the ninth Dedekind number are described. For that, Formal Concept Analysis is used to detect symmetries on lattice intervals. These symmetries could be efficiently computed with the graph isomorphism software nauty.\\

The only thing left to say is that we run the algorithm on Nvidia A100 GPUs. $5311$ GPU hours and $4257682565$ matrix multiplications later, we got the following value for the ninth Dedekind number:

\[286386577668298411128469151667598498812366.\]


\bibliographystyle{splncs04}
\bibliography{FDL9}
\end{document}